\documentclass[12pt]{article}  

    \usepackage{graphicx}
    \usepackage{amsmath}
    \usepackage{amsthm}
    \usepackage{amssymb}
    \usepackage{latexsym}

\theoremstyle{definition}

\newcommand{\BR}{\mbox{$\mathbb{R}$}}
\newcommand{\BN}{\mbox{$\mathbb{N}$}}
\title{\bf Completeness in Probabilistic Metric Spaces}
\author{ Delavar Varasteh $\rm Tafti ^1$, Mahdi $\rm Azhini^2$}
\date{
Department of Mathematics, Science and Research Branch,
 Islamic Azad University, Tehran, Iran.}
\begin{document}
\maketitle
\begin{abstract}
In this paper, we present the Cantor Intersection Theorem and a formulation of
Baire Theorem in complete PM spaces. In addition, the Heine-Borel property
for PM spaces is considered in detail.
\end{abstract}
\section*{1. Introduction}
The idea of probabilistic metric space (briefly PM space) was introduced by
Menger  in [4-6]. Since 1958, B. Schweizer and A. Sklar have been studying
these spaces, and have developed their theory in depth [8-12].
These spaces have also been considered by several other authors [e.g.,
3, 14, 16, 18].  For the historical details, as well as  for the motivation behind the introduction of
 PM spaces, the reader should refer to the book by Schweizer and Sklar [7],
where all the development up to the early 80's are collected.
\footnotetext{(1)Email:delavar.varasteh@gmail.com}
\footnotetext{(2)Corresponding author
:Mahdi Azhini,Email:m.azhini@srbiau.ac.ir}
\paragraph{Definition 1.1.}
A distribution function on $[-\infty ,+\infty]$ is a function $F:[-\infty,
+\infty]\longrightarrow [0,1]$ which is left-continuous on $\BR$,
non-decreasing and $F(-\infty)=0, F(+\infty)=1$. We denote by $\Delta$ the
family  of all distribution functions on $[-\infty,+\infty]$. The order on
$\Delta$ is taken pointwise.
\paragraph{Definition 1.2.}
 The Dirac distribution function, $H_a:[-\infty,+\infty]\longrightarrow [0,1]$
is defined for $a\in[-\infty,+\infty)$ by \\
$H_a(x)=\begin{cases}
0 &\quad  x\in [-\infty,+\infty]\\ 1 &\quad
 x\in (a, +\infty]\end{cases}\;,\\$\\
$H_\infty(x)=\begin{cases} 0 &\quad  x\in [-\infty, +\infty)\\
1 &\quad  x=+\infty\end{cases}.$
\paragraph{Definition 1.3.}
A distance distribution function (d.d.f.) $F:[-\infty,+\infty]\longrightarrow
[0,1]$ is a distribution function with a support contained in $[0, +\infty]$.
The family of all distance distribution functions will be denoted by
$\Delta^+$. Some examples were  considered in [2]. The set $\Delta^+$ can be
partially ordered by the usual pointwise order,
 viz. $F\leq G$ if and only if,
$F(x)\leq G(x)$ for every $x$.
\paragraph{Definition 1.4.}
Let $F$ and $G$ be in $\Delta^+$. Let $h$ be in $(0,1]$, and let $(F,G;h)$
denote the condition $F(x-h)-h\leq G(x)\leq F(x+h)+h$ for all $x$ in
$(-\frac{1}{h},\frac{1}{h})$. The modified l$\acute{e}$vy  distance is the
function $d_L$ defined on $\Delta^+\times \Delta^+$ by
$$d_L(F,G)=\text{inf}\{h|\ \text{both}\  (F,G;h) \text{and}\;
(G,F;h)\text\ {hold}\}.$$\\
The following theorems are proved in [7].
\paragraph\\{Theorem 1.5.}
\begin{itemize}
\item[1)]
The function $d_L$ is a metric on $\Delta^+$.
\item[2)]
If $F$ and $G$ are in $\Delta^+$ and $F\leq G$, then $d_L(G; H_0)\leq d_L(F;
H_0)$.
\end{itemize}
\paragraph{Theorem 1.6.}
 The metric space $(\Delta^+, d_L)$ is compact, and hence complete.
\paragraph{Definition 1.7.}
 A sequence $\{F_n\}$ of d.d.f's converge weakly to a d.d.f $F$ (and we write
$F_n\overset{w}{\longrightarrow} F)$ if and only if the sequence $\{F_n(x)\}$
converges to $F(x)$ at each continuity point $x$ of $F$.

Here we present three theorems that their proofs are in [7].
\paragraph{Theorem 1.8.}
 Let $\{F_n\}$ be a sequence of functions in $\Delta^+$, and let $F$ be in
$\Delta^+$. Then $F_n\overset{w}{\longrightarrow} F$  if and only if
 $d_L(F_n,F)\longrightarrow 0$.
\paragraph{Theorem 1.9.}
Let $h$ be in $(0,1]$. For any $F$ in $\Delta^+$
\begin{align*}
d_L(F,H_0)&=\inf\{h|~(F, H_0; h)\text{\ holds}\}\\
 & = \inf\{h|F(h^+)>1-h\},
\end{align*}
and for any $t>0$, $F(t)>1-t$~ iff~ $d_L(F,H_0)<t$.

\paragraph{Theorem 1.10.}
The supremum of any set of d.d.f's in $\Delta^+$ is in $\Delta^+$.

In order to present the definition of a  probabilistic metric space, we need
the notion of "triangle function" introduced by Serstnev in [13].
\paragraph{Definition 1.11.}
A triangle function $\tau$ is a binary operation on $\Delta^+$ which is
commutative, associative, non-decreasing in each of its variables and has
$H_0$ as identity.

A large class of triangle functions can be constructed through an earlier
concept, which we now introduce.
\paragraph{Definition 1.12.}
A triangular norm (brifly a t-norm) is a binary operation $T$ on the unit
interval $[0,1]$ that is associative, commutative, non-decreasing in each of
its variables and such that $T(x,1)=x$ for every  $x\in [0,1]$. If  $T$ is a
$t$-norm, then its dual $t$-conorm $S:[0,1]^2\rightarrow [0,1]$ is given by
$S(x,y)=1-T(1-x,1-y)$.

It is obvious that a $t$-conorm is commutative, associative, and
non-decreasing operation on $[0,1]$ with unit element $0$.
\paragraph{Example 1.13.}
\begin{itemize}
\item[i)]
Minimm $T_M$ and  maximum $S_M$ given by
$$T_M(x,y)=\text{min}\{x,y\}, \quad S_M(x,y)=\text{max}\{x,y\}.$$
\item[ii)]
product $T_P$ and Probabilistic sum $S_P$ given by
$$T_P(x,y)=xy, \quad S_P(x,y)=x+y-xy.$$
\item[iii)]
Lukasiewicz $t$-norm $T_L$ and Lukasiewicz $t$-conorm $S_L$ given by
$$T_L(x,y)=\text{max}\{x+y-1,0\}, \; S_L(x,y)=\text{min}\{x+y,1\}.$$
\item[iv)]
Weakest $t$-norm (drastic product) $T_D$ and strongest $t$-conorm $S_D$ given
by
\begin{align*}
& T_D(x,y)=\begin{cases}
\text{min}\{x,y\} & \quad
\text{if } \; \text{max}\{x,y\}=1\\
0 & \text{otherwise},
\end{cases}\\
 & S_D(x,y)=\begin{cases}
\text{max}\{x,y\} &\quad  \text{if }\; \text{min}\{x,y\}=0 \\
1 & \text{otherwise}.\end{cases}
\end{align*}
More examples are given in [2].
\end{itemize}
\paragraph{Example 1.14.}
If $T$ is a left-continuous $t$-norm, then the function defined by
$$\tau_T(F,G)(x)=\sup\{T(F(u),G(\nu))|u+\nu=x\},\quad(x\in\BR),$$ is a triangle function (see [7,
section 7.2]).
\paragraph{Example 1.15.}
If $F,G\in \Delta^+$ then we define their  convolution $F*G$ on $[0,
+\infty)$ by\\
 $(F*G)(0)=0,  (F*G)(+\infty)=1$ and
$$(F*G)(x)=\int_{[0,x)}F(x-t)dG(t)\quad \text{for} \; x\in (0,+\infty).$$

The convolution is a triangle function (see[2]).
\section*{2. Probabilistic Metric Space (PM Space)}
\paragraph{Definition 2.1.}
A probabilistic metric space (PM Space) is a triple $(\cal{S},\cal{F},\tau)$
where  $S$ is a nonempty set, ${\cal{F}}:S\times S\rightarrow \Delta^+$ is
given by ${\cal{F}}(p,q)=F_{p,q}\  and\  \tau$ is a triangle function, such that the
following conditions are  satisfied for all $ p,q,r $ in $S$:
\begin{itemize}
\item[a)]
$F_{p,p}=H_0$;
\item[b)]
$F_{p,q}\neq H_0$ for $p\neq q$;
\item[c)]
$F_{p,q}=F_{q,p}$;
\item[d)]
$F_{p,r}\geq \tau(F_{p,q},F_{q,r})$.
\end{itemize}
A very important class of probabilistic metric spaces is given in the
following example,(see[2]).
\paragraph{Example 2.2.}
Let $(S,{\cal{F}},\tau)$ be PM space and $\tau=\tau_T$, where
$$\tau_T(F,G)(x)=\sup\{T(F(u),G(\nu))|u+\nu=x\}, \quad (x\in\BR),$$
for a $t$-norm $T$. Then $(S,{\cal{F}},\tau)$ is called Menger space.\\
As a very  special case of a Menger space we obtain the classical metric space.
\paragraph{Example 2.3.}
 We can prove that $(M,{\cal{F}},\tau_T)$ , for ${\cal{F}}_{(p,q)}=H_{d(p,q)}$, for all $ p,q $ in $ M$ and any $t$-norm $T$ is a Menger space  if and only if $(M,d)$ is a classical metric space,(see [2]).
\paragraph{Example 2.4.}
 A PM space $(S,{\cal{F}},\tau)$ for which $\tau$ is a convolution is called
Wald space.\\Many different topological structures may be defined on a PM space. The one
that has received the most attention to date is the strong topology.
\paragraph{Definition 2.5.}
Let $(S,\cal{F},\tau)$ be a PM space. The strong topology is introduced  by a strong neighbourhood system ${\cal{N}}=\bigcup_{p\in S}{\cal{N}}_p$, where
${\cal{N}}_p=\{N_p(t)|t>0\}$ and {$N_p(t)=\{q\in S|F_{p,q}(t)>1-t\}$}for $t>0$
and $p\in S$.\\ Applying Theorem (1.9), we see that $q$ is in $N_p(t)$ if and
only if $d_L(F_{p,q},H_0)<t$, whence $N_p(t)=\{q\in S|d_L(F_{p,q}, H_0)<t\}$.\\
The $(\epsilon , \lambda)$-topology on
$(S,\cal{F},\tau)$ which is introduced by a family of  neighbourhood\\
$\{N_p(\epsilon, \lambda)\}_{p\in S, \epsilon>0, \lambda\in (0,1)}$, where
$N_p(\epsilon ,\lambda)=\{q\in S|F_{p,q}(\epsilon)>1-\lambda\}$ is of special interest.
Since $N_p(t,t)=N_p(t)$ for $t>0$ and $N_p(\text{min}\{\epsilon ,
\lambda\})\subseteq N_p(\epsilon ,\lambda)$ for every  $\epsilon>0$,
$\lambda\in (0,1)$; thus  the strong neighbourhood system and the  $(\epsilon
, \lambda)$ -neighbourhood system define the same topology. (see [7, section
4.3]).
\paragraph{Theorem 2.6([7]).}
 Let $(S,\cal{F}, \tau)$ be a PM space. If $\tau$ is continuous, then the strong
neighbourhood system $\cal{N}$ satisfies (1), (2).
\begin{itemize}
\item[(1)]
If $V$ is a strong  neighbourhood of $p\in S$, and $q$ is in $V$, then there
is a strong neighbourhood $W$ of $q$ such that $W\subseteq V$.
\item[(2)]
If $p\neq q$, then there is a $V$ in ${\cal{N}}_p$ and  a $W$ in ${\cal{N}}_q$
such that  $V\cap W=\phi$ and thus the strong  neighbourhood system $\cal{N}$
determines a Hausdorff topology for $S$.
\end{itemize}
\paragraph{Definition 2.7.}
 Let $(S,\cal{F}, \tau)$ be a  PM space.
 A sequence $\{p_n\}_{n\in\BN}$ in $S$ converges to $p\in S$
 in $(\epsilon , \lambda)$-topology if for every $\epsilon >0$ and
$\lambda\in (0,1)$, there exists~
 $n_0(\epsilon , \lambda)\in \BN$~ such that \\
$F_{p_n,p}(\epsilon)>1-\lambda$ for every $n\geq n_0(\epsilon ,\lambda)$.
\paragraph{Definition 2.8.}
 Let $(S,\cal{F}, \tau)$ be a  PM space.
 A sequence $\{p_n\}_{n\in\BN}$ in $S$  is a Cauchy sequence if
 for every $\epsilon >0$ and
$\lambda\in (0,1)$, there exists $n_0(\epsilon , \lambda)\in \BN$ such that
$F_{p_n,p_m}(\epsilon)>1-\lambda$ for every $n,m\geq n_0(\epsilon ,\lambda)$,
and $S$ is  complete if every Cauchy sequence in $S$ converges to a point in
 $S$.
\paragraph{Theorem 2.9([7]).}
 Let $(S,\cal{F}, \tau)$ be a  PM space.
 A sequence $\{p_n\}$ in $S$  converges to  a point $p$ in $S$ if and only if
$d_L(F_{p_n,p},H_0)\rightarrow 0$, as $n \rightarrow \infty$.Similarly, a sequence $\{p_n\}$ in $S$ is a
 Cauchy sequence if  and only if $d_L(F_{p_m,p_n},H_0)\rightarrow 0$,\ as $m,n\rightarrow\infty$.
\paragraph{Definition 2.10.}
 Let $(S,\cal{F}, \tau)$ be a  PM space.
The probabilistic diameter $D_A$ of a nonempty subset $A$ of $S$ is the
function $D:[-\infty,+\infty]\rightarrow \BR$  defined by $D_A(+\infty)=1$ and
 for every $x>0$
$$D_A(x)=\text{sup}_{t<x}\,[{\inf}_{p,q\in A}F_{p,q}(x)].$$
\paragraph{Definition 2.11.}
 Let $(S,\cal{F}, \tau)$ be a  PM space and $A$ a nonempty subset of $S$.
\begin{itemize}
\item[1)]
A is said to be probabilistic bounded if $\sup_{x>0} D_A(x)=1.$
\item[2)]
A is said to be probabilistic semi-bounded if $0<\sup_{x>0}D_A(x)<1$.
\item[3)]
A is said to be probabilistic unbounded if $D_A(x)\equiv 0.$
\end{itemize}
\paragraph{Definition 2.12.}
 Let $(S,\cal{F}, \tau)$ be a  PM space and $A$ a nonempty subset of $S$.
If $\{V_\alpha\}_{\alpha\in I}$ is a family of sets such that  $A\subset
\cup_{\alpha\in I} V_\alpha,\{V_\alpha\}_{\alpha\in I}$ is called  a cover of
$A$, and $A$ is said to be covered  by $V_\alpha ^,  s$.
A is said to be probabilistic  totally bounded if, for every $\epsilon>0$, $A$
can be  covered by finitely  many neighbourhoods $N_p(\epsilon)$,  where $p\in
A$ or equivalently if for each $\epsilon >0$ and for each $p\in A$, there
exists $p_1,p_2,\cdots , p_n\in A$ such that $F_{p_i,p}(\epsilon )=1$ for
every $i\in\{1,2,\cdots ,n\}$ (see[1],[15]).

We now establish the properties of the probabilistic diameter.
\paragraph{Theorem 2.13([7]).}
 Let $(S,\cal{F}, \tau)$ be a  PM space  under a continuous triangle function
$\tau$ and $A$ a nonempty  subset of $S$. The probabilistic diameter $D_A$ has
the following properties:
\begin{itemize}
\item[1)]
The function $D_A$ is a distance distribution function.
\item[2)]
$D_A=H_0$ iff $A$ is  a singleton set.
\item[3)]
 if $A\subseteq B$, then $D_A\geq D_B$.
\item[4)]
For any $p,q$ in $A$, $F_{p,q}\geq D_A$.
\item[5)]
 If $A=\{p,q\}$, then $D_A=F_{p,q}$.
\item[6)]
If $A\cap B$ is nonempty, then $D_{A\cup B} \geq \tau(D_A,D_B)$.
\item[7)]
$D_A=D_{\bar{A}}$, where $\bar{A}$ denotes the closure of $A$ in the
$(\epsilon ,\lambda)$- topology on $S$.
\end{itemize}
\paragraph{Theorem 2.14.}
In any PM space every probabilistic totally bounded set is probabilistic bounded.
\paragraph{Proof.}
 Let $(S,\cal{F}, \tau)$ be a PM space. Let $A\subseteq S$ be  a
probabilistic totally bounded set. Hence for every $p\in A$ and $\epsilon >0$, there exists $p_1,\cdots , p_n\in
A$ such that $F_{p_i,p}(\epsilon)=1$ for every  $i\in\{1,2,\cdots , n\}$.
According to Definition 2.1 (d),  we have $1=\tau(F_{p,p_i}(\epsilon),
F_{p_i,q}(\epsilon))\leq F_{p,q}(\epsilon)\leq 1$ for every $\epsilon>0$,
$p,q\in A$ and $i\in\{1,2,\cdots, n\}$. Hence $F_{p,q}(\epsilon)=1$ for each
$p,q\in A$ and  $\epsilon >0$. Therefore  $\sup_{\epsilon>0}\inf_{p,q\in A}
F_{p,q}(\epsilon)=1$, that is $A$ is probabilistic bounded. \hskip 3.8cm $\Box$
\paragraph{Theorem 2.15 ([17])}
 Let $(S,\cal{F}, \tau)$ be a  PM space  and $\tau$ be continuous. Let
$\{p_n\}_{n\in\BN}$ be a sequence of $S$. Also suppose that:
\begin{itemize}
\item[a)]
$\{p_n\}_{n\in\BN}$ is a Cauchy sequence,
\item[b)]
there exists a  subsequence $\{p_{i_n}\}_{n\in\BN}$
 of $\{p_n\}_{n\in\BN}$ that
converges  to some $p_0\in S$.
\end{itemize}
Then $\{p_n\}_{n\in\BN}$ converges to $p_0$.
\paragraph{Theorem 2.16 (The Cantor Intersection Theorem).}
 A probabilistic metric space $(S,\cal{F},\tau)$ under a continuous traingle
function $\tau$ is complete if and only if for every  nested sequence
$\{S_n\}_{n\in\BN}$  of nonempty closed subsets  of $S$,  such that
$D_{S_n}\rightarrow H_0$ (pointwise) as $n\rightarrow \infty$, the
intersection $\bigcap_{n=1}^\infty S_n$ consists of exactly one point.
\paragraph{Proof.}
First suppose that $S$ is complete and $\{S_n\}_{n\in\BN}$ is a nested
sequence of nonempty closed subsets of $S$ such that $D_{S_n}\rightarrow
H_0$ pointwise as $n\rightarrow \infty$.
 So there is a positive integer  $n_0(\epsilon ,\lambda)$ such that
$D_{S_n}(\epsilon)>1-\lambda$ whenever $n\geq n_0(\epsilon ,\lambda)$, for all
$\epsilon >0$ and $\lambda\in (0,1)$. Let $m>n\geq n_0(\epsilon,\lambda)$,
then   $p_m\in S_m\subseteq S_n$ and $p_n\in S_n$. Thus by Theorem 2.13(4),
\;$F_{p_n,p_m}(\epsilon )\geq D_{S_n}(\epsilon)>1-\lambda$ and since $\epsilon ,
\lambda$ are arbitrary, hence $\{S_n\}$ is a Cauchy sequence in
$(S,\cal{F},\tau)$. The completeness of $(S,\cal{F},\tau)$ guarantees that
there is an element $p_0$ in $S$ such that $p_n\rightarrow p_0$ as
$n\rightarrow \infty$. Since $p_m\in S_n$ whenever  $m>n$ and $S_n$ is closed,
it follows  that $p_0\in S_n$ for every positive integer $n$. To show the
uniqueness
 of $p_0$, suppose that $q_0$ belonging  to $\bigcap_{n=1}^\infty
S_n$. Apply Theorem 2.13 (4) again to obtain $F_{p_0,q_0}\geq D_{s_n}$.
Now $D_{s_n}\rightarrow H_0$ as $n\rightarrow \infty$, therefore
$F_{p_0,q_0}(\epsilon)=1$, for every $\epsilon>0$ i.e, $p_0=q_0$.

 Conversely,
suppose every nested sequence of nonempty of closed subsets of $S$ with  the
probabilistic diameter tending to  $H_0$, have
an intersection consists of exactly one point. Suppose
$\{p_n\}$ be a Cauchy sequence in $(S,\cal{F},\tau)$. Put
$S_n=\overline{\{p_n,p_{n+1},\cdots\}}$, where $n=1,2,\cdots$. Every $S_n$ is
nonempty and closed subsets  in $S$. $\{S_n\}_{n\in\BN}$ is a nested sequence
of subsets of $S$, therefore  by Theorem 2.13 (3), $\{D_{s_n}\}_{n\in\BN}$ is
an increasing sequence and hence $D_{s_n}\rightarrow H_0$ as $n\rightarrow
\infty$. So there  is $p_0\in S$ such that $\bigcap_{n=1}^\infty S_n=\{p_0\}$.
 Therefore there exists a subsequence  of $\{p_n\}$ such that converges to
$p_0\in S$. Hence by Theorem 2.15, $\{p_n\}$ converges to $p_0\in S$.\hskip 4 cm $\Box$
\paragraph{Theorem 2.17 (Baire Theorem).}
Let $(S,\cal{F},\tau)$ be a complete PM Space  under a continuous triangle
function $\tau$. If $\{G_n\}$ is a sequence of dense and open subset of $S$,
then $\bigcap_{n=1}^\infty G_n$ is not empty. (In fact it is dense in $S$).
\paragraph{Proof.}
$G_1$ is dense in $S$,
 therefore $G_1$ is nonempty. We can choose $p_1\in G$,
since $G_1$ is open in $S$, there is $r>0$ such that $N_{p_1}(r)\subseteq
 G_1$.
First we prove $\overline{N_{p_1}}(r_1)\subseteq N_{p_1}(r)$.
Let $0<r_1<\frac{r}{2}$. Let $p\in \overline{N_{p_1}}(r_1)$, then there exists a
sequence $\{p_n\}$ in $N_{p_1}(r_1)$ such that converges to $p$.
Since $p_n\longrightarrow p$, then for all $\alpha>0$, there exist integer
$n_0(\epsilon)$ such that $d_L(F_{p_n,p}, H_0)<\alpha$  for every $n\geq
n_0(\epsilon)$. By Theorem 1.6, the metric space $(\Delta^+,d_L)$ is compact
and hence $\tau$ is uniformly continuous on $\Delta^+\times \Delta^+$ and
since  $\tau$ has $H_0$ as identity we have $\tau(H_0,F_{p_n,p_1})=F_{p_n,p_1}$.
 Thus for all $\epsilon>0$ , there is $\delta>0$ such that
$d_L(\tau(F_{p_n,p},F_{p_n,p_1}),F_{p_n,p_1})<\epsilon$ whenever
$d_L(F_{p_n,p},H_0)<\delta$.In the above argument, put $\alpha=\delta$. On the
other hand, according
to Theorem 1.5 (2) we have
$$d_L(F_{p,p_1}, H_0)\leq d_L(\tau(F_{p_n,p},F_{p_n,p_1}),H_0).$$
By Theorem 1.5 (1), $d_L$ is a metric on  $\Delta^+$, therefore  we  obtain
\begin{align*}
d_L(F_{p_1,p},H_0) & \leq d_L(\tau(F_{p_n,p},F_{p_n,p_1}),H_0)\\
& \leq  d_L(\tau(F_{p_n,p},F_{p_n,p_1}),F_{p_n,p_1})+d_L(F_{p_n,p_1},H_0)\\
& \leq \epsilon +d_L(F_{p_n,p_1},H_0).
\end{align*}
Since $p_n\in N_{p_1}(r_1)$, therefore $d_L(F_{p_n,p_1},H_0)<r_1$. Now we
choose $\epsilon =r_1-d_L(F_{p_n,p_1},H_0)$, hence we conclude
\begin{align*}
d_L(F_{p,p_1},H_0) & < \epsilon +d_L(F_{p_n,p_1},H_0)\\
& = r_1-d_L(F_{p_n,p_1},H_0)+d_L(F_{p_n,p_1},H_0)\\
& = r_1<r.
\end{align*}
 Thus  $p\in N_{p_1}(r)$ i.e. $\overline{N_{p_1}}(r_1)\subseteq N_{p_1}(r)$.As
$G_2$ is dense in $S$ by Theorem 2.6 and $N_{p_1}(r_1)$ is open in $S$, hence
$G_2\cap N_{p_1}(r_1)$ is nonempty, therefore we can choose $p_2\in G_2\cap
N_{p_1}(r_1)$.
By  Theorem 2.6, intersection of two open sets is open, thus there  is
$0<r_2<\frac{r_1}{2}$ such that $\overline{N_{p_2}}(r_2)\subseteq  G_2\cap
N_{p_1}(r_1)$. By induction,  we choose $p_n\in S$ and $r_n>0$ as  follows: with
$p_i, r_i$ if $i<n$, we see that $G_n\cap N_{p_{n-1}}(r_{n-1})$ is nonempty
and open, therefore we can  choose $p_n,r_n$ so that $r_n<\frac{r}{2^n}$ and
$\overline{N_{p_n}}(r_n)\subseteq G_n\cap N_{p_{n-1}}(r_{n-1})$.

Every $\overline{N_{p_n}}(r_n)$ is nonempty and closed strongly in $S$.
$\{\overline{N_{p_n}}(r_n)\}$ is a nested sequence of subsets of $S$. Also
$D_{\overline{N_{p_n}}(r_n)}\longrightarrow H_0$  as $n\rightarrow \infty$.
Because, by the definition of neighbourhood , we have
$N_{p_n}(r_n)=\{q\in S|d_L(F_{p_n,q}, H_0)<r_n\}$.Since $r_n<\frac{1}{2^n}r$, for all $n\geq 1$,
we see that $d_L(F_{p_n,q},H_0)\rightarrow  0$ as $n\rightarrow \infty$.
Therefore for each $x>0$, we have $F_{p_n,q}(x)\rightarrow 1$ as $n\rightarrow \infty$
i.e., $p_n=q$ for all $n\geq 1$. Thus $N_{p_n}(r_n)$ is singleton.
By Theorem 2.13 (7), we have $D_{\overline{N}_{p_n}(r_n)}=D_{N_{p_n}(r_n)}$,
therefore according to Theorem 2.13 (2), $D_{\overline{N}_{p_n}(r_n)}
 \rightarrow
H_0$ as $n\rightarrow \infty$. Now by  Theorem 2.16 (The Cantor Intersection Theorem),$\bigcap_{n=1}^\infty{\overline{N}_{p_n}(r_n)}$ is not empty and since for all $n\geq 1$,we have $\bigcap_{n=1}^\infty{\overline{N}_{p_n}(r_n)}\subseteq \bigcap_{n=1}^\infty G_n$,therefore
$\bigcap_{n=1}^\infty G_n$ is not empty. In fact $\bigcap_{n=1}^\infty G_n$ is
dense in $S$.As $G_1$ is  dense in $S$, therefore  $N_p(r)\cap G_1$ is not
empty, for every $p\in S$  and every $r>0$. We put $A_1=N_p(r)\cap G_1$, hence
$A_1$ is open in $S$, thus it contains a  neighbourhood $N_{p_1}(t)$ for every
$p_1\in A_1$.According to the preceding argument, there is $r_1>0$  such that
$\overline{N_{p_1}}(r_1)\subseteq N_{p_1}(t)\subseteq  A$.
Now we set $A_n=G_n$ for $n=2,3,\cdots $ and $E_n=\overline{N_{p_n}}(r_n)$ so, we
have  $\bigcap_{n=1}^\infty E_n\subseteq \bigcap_{n=1}^\infty A_n$, therefore
$\bigcap_{n=1}^\infty A_n$ is not empty. On the other hand,
$\bigcap_{n=1}^\infty A_n=A_1\cap(\bigcap_{n=2}^\infty A_n)$ , thus
$(N_p(r)\cap G_1)\cap (\bigcap_{n=2}^\infty G_n)\neq \phi$, therefore
$N_p(r)\cap (\bigcap_{n=1}^\infty G_n)\neq \phi$. Because $N_p(r)$ is an
arbitrary open set in $S$,hence $\bigcap_{n=1}^\infty G_n$ is dense in $S$.\hskip 16.29 cm $\Box$

\paragraph{Theorem 2.18.}
Suppose $(S,\cal{F},\tau)$ is a PM  space  and $\tau$  is continuous. If $E$
is a subset of the PM space $(S,\cal{F},\tau)$, the following are equivalent:
\begin{itemize}
\item[a)]
$E$ is complete and probabilistic totally bounded.
\item[b)]
(The Bolzano-Weierstrass Property). Every sequence in $E$  has a subsequence
converges to a point of $E$.
\item[c)]
 (The Heine-Borel Property). If $\{V_\alpha\}_{\alpha\in A}$ is a cover of $E$
by open sets, then there is a finite set $F\subset A$ such that
$\{V_\alpha\}_{\alpha\in F}$ covers  $E$.
\end{itemize}
\paragraph{Proof.}
We first show that (a) and (b) are equivalent. Then (a) implies(c),  and finally (c) implies (b).

(a) implies (b): Suppose that (a)
holds and  $\{x_n\}$ is a sequence in $E$. Since $E$ is probabilistic totally
bounded, hence it can be  covered by finitely many neighbourhoods
$N_{p_i}(2^{-1})$, with $p_i\in E$ and $i\in \{1,2,\cdots ,n\}$. At least one
of them must  contain $x_n$ for infinitely  many $n$: say, $x_n\in
N_{p_1}(2^{-2})$ for $n\in N_1$. $E\cap N_{p_1}(2^{-1})$ can be  covered by
finitely many neighbourhoods $N_{p_i}(2^{-2})$, with $p_i\in E$ and
$i\in\{1,2,\cdots ,n\}$ and  at least one of them must contain $x_n$ for
infinitely  many $n\in N_1:$ say, $x_n\in N_{p_2}(2^{-2})$ for $n\in N_2$.
Continuing inductively, we obtain a sequence of neighbourhoods
$N_{p_j}(2^{-j})$, with $p_j\in E$ and a decreasing sequence of subset $N_j$
of $\BN$  such that $x_n\in N_{p_j}(2^{-j})$ for $n\in N_j$. Pick $n_1\in
N_1$, $n_2\in N_2,\cdots $ such that $n_1<n_2<\cdots $, then $\{x_{n_j}\}$ is a
Cauchy sequence. Because, by Theorem 1.5 (2) we have
$$d_L(F_{x_{n_j}, x_{n_k}}, H_0)\leq  d_L (\tau(F_{x_{n_j},x_j},
F_{{x_j},x_{n_k}}), H_0).$$
By Theorem 1.5 (1) $d_L$ is a metric on $\Delta^+$, therefore we obtain
\begin{align*}
d_L(F_{x_{n_j},x_{n_k}},H_0) & \leq
d_L(\tau(F_{x_{n_j},x_j},F_{x_j,x_{n_k}}),H_0)\\
& \leq d_L (\tau(F_{x_{n_j},x_j}, F_{x_j,x_{n_k}}),F_{x_j,x_{n_k}})+d_L
(F_{x_j,x_{n_k}},H_0).
\end{align*}
 By Theorem 1.6, the metric space $(\Delta^+,d_L)$ is compact and hence $\tau$
is uniformly continuous on $\Delta^+\times \Delta^+$ and  since $\tau$ has
$H_0$ as identity i.e., $\tau(H_0,F_{x_{n_j}, x_j})=F_{x_{n_j},x_j}$, thus for all $\epsilon>0$
there is $\delta>0$ such that
$d_L(\tau(F_{x_{n_j},x_j},F_{x_j,x_{n_k}}),F_{x_{n_j},x_j})<\epsilon$ wherever
$d_L (F_{x_j,x_{n_k}},H_0)<\delta$. Since $E$ is
probabilistic totally bounded, hence $d_L(F_{x_j,x_{n_k}}, H_0)<\epsilon^\prime$, for
all $\epsilon^\prime >0$. We put $\epsilon^\prime =\delta$, then $\{x_{n_j}\}$ is Cauchy
sequence, for $d_L(F_{x_{n_j}, x_{n_k}}, H_0)\leq 2^{1-j}$, if $k>j$, and
since $E$ is complete, it has a limit in $E$.

(b) implies (a):  We show that if any conditions in (a) fails, then so  does
(b). If $E$ is not complete, there is a Cauchy sequence $\{x_n\}$ in $E$ with
 no limit in $E$ so,none of subsequences of $\{x_n\}$ can converge in $E$, otherwise
the  whole sequence would converges the same limit. Since $\{x_n\}_{n\in\BN}$
is a Cauchy  sequence, for all $\epsilon>0$ there is $N_1\in\BN$  such that $d_L(F_{x_n,x_m},
H_0)<\frac{\epsilon}{2}$ whenever $n\geq N_1$, $m\geq N_1$
. Now if there exists a subsequence $\{x_{i_n}\}$ of $\{x_n\}$  that
converges to some $x_0\in E$ then,for all $\epsilon>0$ there is $N_2\in\BN$ such that
$d_L(F_{x_{i_n},x_0},H_0)<\frac{\epsilon}{2}$ whenever $n\geq N_2$
. By uniform continuity of $\tau$ ,we have
\begin{align*}
d_L(F_{x_n,x_0},H_0)&\leq d_L(\tau(F_{x_n,x_{i_n}},F_{x_{i_n},x_0}),H_0)\\
& \leq
d_L(\tau(F_{x_n,x_{i_n}},F_{x_{i_n},x_0}),F_{x_n,x_{i_n}})+d
_L(F_{x_n,x_{i_n}},H_0)\\
& <\frac{\epsilon }{2}+\frac{\epsilon}{2}=\epsilon.
\end{align*}
  Therefore $\{x_n\}$ converges to $x_0$, a contradiction. On the other hand, if $E$ is not
probabilistic totally bounded, let $\epsilon>0$ be such that $E$ cannot  be
covered by finitely many neighbourhoods $N_{x_i}(\epsilon)$, with $x_i\in E$
and $i\in \{1,2,\cdots ,n\}$. Choose $x_n\in E$ inductively as follows. Begin
with any $x_1\in E$, and having chosen $x_1,\cdots , x_n$, pick $x_{n+1}\in
E\backslash \bigcup_{i=1}^nN_{x_i}(\epsilon)$ for all $n\in\BN$. This
implies, $d_L(F_{x_m,x_n},H_0)\geq \epsilon $ if $m\neq n$. So $\{x_n\}$ has
no Cauchy subsequence which contradicts with Bolzano-Wierstrass property.\\
(b) implies(c): It suffices to show that if (b) holds and
$\{V_\alpha\}_{\alpha\in A}$ is a cover of $E$ by open sets, there exists
$\epsilon> 0$ such that every neighbourhood $N_p(\epsilon)$, with $p\in E$
that intersects $E$,is contained in some $V_\alpha$, for $E$ can be covered by
finitely many such neighbourhoods by (a). Suppose to the contrary, for every
$n\in\BN$, there is $x_n\in E$ such that $N_{x_n}(2^{-n})$ is not contained in
any $V_\alpha$, with $\alpha\in A$. By the Bolzano-Weierstrass property, there
exists a subsequence $\{x_{i_n}\}$ that converges to some $x\in E$. Since
$\{V_\alpha\}_{\alpha\in A}$ is a covering of $E$, there is $\alpha\in A$ such
that  $x\in V_\alpha$.
 Also $V_\alpha$ is open therefore there exists $\epsilon
>0$ such that $N_x(\epsilon)\subset V_\alpha$. Now we  prove
$N_{x_{i_n}}(2^{-n})\subseteq V_\alpha$. If $y\in N_{x_{i_n}}(2^{-n})$, then
$d_L(F_{x_{i_n},y}, H_0)<2^{-n}$. By the uniform continuity of $\tau$, for all $\epsilon >0$ there
exists  $\delta>0$ such that
$d_L(\tau(F_{x,x_{i_n}},F_{x_{i_n},y}),F_{x_{i_n},y})<\epsilon$ whenever
$d_L(F_{x_{i_n},x},H_0)<\delta$. Now let $n\in \BN$ be
such that $d_L(F_{x_{i_n},x},H_0)<\delta$ and $2^{-n}<\frac{\epsilon}{2}$.
Therefore  we have
\begin{align*}
d_L(F_{x,y},H_0) & \leq d_L (\tau(F_{x,x_{i_n}},F_{x_{i_n},y}),H_0)\\
& \leq
d_L(\tau(F_{x,x_{i_n}},F_{x_{i_n},y}),F_{x_{i_n},y})+d_L(F_{x_{i_n},y},H_0)\\
& < \frac{\epsilon }{2}+2^{-n}<\epsilon.
\end{align*}
 Hence $N_{x_{i_n}}(2^{-n})\subseteq V_\alpha$, which is a contradiction.

(c) implies (b): If $\{x_n\}$ is a sequence in $E$ with no convergent
subsequence, for each $x\in E$ there is $\epsilon(x)>0$ such that $N_{x}{(\epsilon(x))}$
 contains $x_n$ for only finitely many $n$. Otherwise some
subsequence of $\{x_n\}$  converges to $x$. Then $\{N_{x}{(\epsilon(x))}\}_{x\in E}$
 is a cover of $E$ by open sets with no finite subcover.\hskip15 cm $\Box$


\begin{thebibliography}{99}
\bibitem{}
Engelking, R;  General Topology, Warszawa (1977).
\bibitem{}
 O. Had$\check{z}$i$\acute{c}$
 and  E. Pap, Fixed  point theory in probabilistic
metric spaces, Institue of Mathematics, University of Novi Sad, Yugoslavia
(2001).
\bibitem{}
 V. Istratescu and I. Vaduva, Products of statistical metric spaces
(Romanian), Acad.R.P.Romine Stud.Cere.Mat.12(1961),567-574.
\bibitem{}
K. Menger, Statistical metrics, Proc, Nat, Acad, of sci. U.S.A. 28 (1942),
535-537.
\bibitem{}
K. Menger, Probabilitstic geomety, Pro,Nat, Acad,  of sci, U.S.A. 37(1951),
226-229.
\bibitem{}
K. Menger, Gemetric generale (chap.VII), Memorial des Sciences Mathematiques, No, 124, Paris 1954.
\bibitem{}
B. Schweizer and A. Sklar, probabilistic metric spaces, North-Holland, New York(1983).
\bibitem{}
B.  Schweizer and A. Sklar, Espaces
 metriques aleatoires, C.R.Acad. Sci., Paris 247(1958),2092-2094.
\bibitem{}
B. Schweizer and A. Sklar, Statistical metric spaces, Pacific J. Math. 10 (
1960), 314-334.
\bibitem{}
B. Schweizer and A. Sklar, Associative functions and statistical triangle
inequalitles, Publicationes Mathematicae, Debrecen 8(1961), 169-186.
\bibitem{}
B. Schweizer and A. Sklar, Triangle inequalities in a class of statistical
metric spaces, J. London Math. Soc. 38(1963), 401-406.
\bibitem{}
B. Schweizer and A. Sklar and E. Thorp, The metrization of statistical metric
Spaces, Pacific J. Math. 10(1960), 673-675.
\bibitem{}
A. N. $\check{S}$erstnev, Random normed spaces: Problems of completeness,
Kazan. Gos.  Univ. U$\check{c}$en. Zap. 122(1962), 3-20.
\bibitem{}
A. N. $\check{S}$erstnev, On the concept of random normed spaces, (Russian),
Dokl. Akad, Nauk SSSR 149(1963) 280-283.
\bibitem{}
Zhang Shi-Sheng, Basic theory and application of probabilistic metric space
(I), Applied Mathematics and Mechanics, 9,2 (1988); 123-133.
\bibitem{}
E. Thorp, Generalized topologies for statistical metric spaces, Fundamenta
Mathematicae 51 (1962), 9-21.
\bibitem{}
S. M. Vaezpour and M. Shams, Topological properties of probabilistic metric
space Int. Journal of Math. Analysis, Vol. 1, 2007, no. 20, 957-964.
\bibitem{}
A. Wald, On a statistical generalization of metric spaces, Proc, Nat. Acad. of
Sci U. S. A. 29 (1943), 196-197.
\end{thebibliography}
\end{document}